\documentclass[11pt]{amsart}
\input{epsf.tex}
\usepackage{euscript,color,latexsym,amsfonts,amssymb,epsfig,verbatim}
\usepackage{amsmath,amsthm,latexsym,graphics,textcomp}
\usepackage{url}
\usepackage{graphicx}

\input xy
\xyoption{all}

\theoremstyle{plain}
\newtheorem{theorem}{Theorem}[section]
\newtheorem{lemma}[theorem]{Lemma}

\newtheorem{remark}[theorem]{Remark}

\newtheorem*{claim}{Claim}

\def\Mod{\mbox{\rm{Mod}}}

\def\SL{\mbox{\rm{SL}}}
\def\GL{\mbox{\rm{GL}}}

\def\PSL{\mbox{\rm{PSL}}}

\def\Aff{\mbox{\rm{Aff}}}

\def\PML{{\mathbb P}{\mathcal {ML}}}

\def\GL{{\mathcal {GL}}}

\def\Cent{\mbox{\rm{C}}}

\begin{document}

\title{Accidental parabolics in mapping class groups}

\author{Christopher J Leininger}\thanks{Research supported by NSF DMS 06-03881}

%\address{Department of Mathematics \\ Columbia University \\ 2990 Broadway MC 4448 \\
%New York, NY 10027-6902}

%\thanks{The author was supported in part by an NSF postdoctoral fellowship.}

%\email{clein@math.columbia.edu}

%\keywords{Veech groups, mapping class group, parabolic}

%\subjclass[2000]{Primary: ????; Secondary: ????}

\maketitle

\begin{abstract}
In this note we discuss the behavior of the Gromov boundaries and limit sets for the surface subgroups of the mapping class group with accidental parabolics constructed by the author and A. Reid in \cite{leiningerreid}.  Specifically, we show that generically there are no Cannon--Thurston maps from the Gromov boundary to Thurston's boundary of Teichm\"uller space.
\end{abstract}

%%%%%%%%%%%%%%%%%%%%%%%%%%%%%%%%%%%%%%%%%%%
%%%%%%%%%%%%%%%%%%%%%%%%%%%%%%%%%%%%%%%%%%%
\section{Introduction} \label{introsect}
%%%%%%%%%%%%%%%%%%%%%%%%%%%%%%%%%%%%%%%%%%%
%%%%%%%%%%%%%%%%%%%%%%%%%%%%%%%%%%%%%%%%%%%
%%%%%%%%%%%%%%%%%%%%%%%%%%%%%%%%%%%%%%%%%%%

Let $S = S_g$ be a closed genus $g$ surface ($g \geq 2$) and $\Mod(S)$ its mapping class group. In
\cite{leiningerreid}, we consider a family of subgroups $G(\omega_g) < \Mod(S_g)$ depending on a certain holomorphic
abelian differential $\omega_g$ constructed by Veech \cite{veech1,veech2} (we recall the relevant geometry of these
groups below). The group $G(\omega_g)$ is free on $2g$--generators but is naturally isomorphic to the fundamental group
of a surface of genus $g$ with one puncture.  Indeed, $G(\omega_g)$ stabilizes a totally geodesic hyperbolic plane
$\mathbb H_{\omega_g}$ in Teichm\"uller space with quotient a one--cusped hyperbolic surface.  We write $G_0(\omega_g)
< G(\omega_g)$ for a representative of the unique conjugacy class of maximal parabolic subgroups.

For each $g$, we can abstractly double $G(\omega_g)$ over $G_0(\omega_g)$ to obtain $\mathcal G_{2g}$ the fundamental group of a closed
surface of genus $2g$
$$\mathcal G_{2g} = G(\omega_g) *_{G_0(\omega_g)} G(\omega_g).$$

We can attempt to carry out this doubling within $\Mod(S)$ by conjugating one of the factors of $\mathcal G_{2g}$ by an
element of $\Cent(G_0(\omega_g))$, the centralizer of $G_0(\omega_g)$ in $\Mod(S)$. More precisely, for every $h \in
\Cent(G_0(\omega_g))$ and $n \in \mathbb Z$ we obtain a homomorphism from $\mathcal G_{2g}$ to the subgroup of $\Mod(S)$
generated by $G(\omega_g)$ and its conjugate by $h^n$
$$\mathcal G_{2g} \to \mathcal G(\omega_g,h,n) = \langle G(\omega_g),h^n G(\omega_g) h^{-n} \rangle < \Mod(S).$$
This homomorphism is the canonical one extending the obvious isomorphisms of the first and second factors to $G(\omega_g)$ and $h^n G(\omega_g) h^{-n}$, respectively.

For arbitrary $h \in \Cent(G_0(\omega_g))$ and $n \in \mathbb Z$ this homomorphism need not be injective. This setup
however is reminiscent of the Maskit Combination Theorem in hyperbolic geometry \cite{maskit}. In \cite{leiningerreid},
we show that under mild hypothesis, one is able to make a similar conclusion to that of Maskit's Theorem.  To state the result precisely, we recall that $G_0(\omega_g)$ is (virtually) generated by a multitwist in a
multicurve $A$  (see Section \ref{examplesection}).

\begin{theorem} [{\bf Leininger--Reid} \cite{leiningerreid}] \label{lrthrm}
Suppose $h \in \Cent(G_0(\omega_g))$ is pseudo-Anosov on $S - A$. Then
$$\mathcal G_{2g} \to \mathcal G(\omega_g,h,n) < \Mod(S_g)$$
is an isomorphism for all sufficiently large $n$.  Moreover, every element of $\mathcal G(\omega_g,h,n)$ is
pseudo-Anosov except those conjugate into $G_0(\omega_g)$.
\end{theorem}

The proof of Theorem \ref{lrthrm} in \cite{leiningerreid} mimics the proof of Maskit's Theorem, at least in principle, and
the groups $\mathcal G(\omega_g,h,n)$ are analogous to Kleinian surface groups with an {\em accidental parabolic}.
One may wonder the extent to which these groups behave like their Kleinian counterparts. For
example, one can ask if they satisfy some form of geometric finiteness---see Problem 6.1 of \cite{mosherproblems}. One
ingredient that seems desirable for such a notion is a description of the ideal boundary behavior---compare Floyd
\cite{floyd} for the Kleinian setting. In \cite{leiningerreid} Question 10.1, we ask whether or not there is a
\textbf{Cannon--Thurston map} for these groups
$$\partial \mathcal G(\omega_g,h,n) \to \PML(S_g)$$
By this, we mean a $\mathcal G(\omega_g,h,n)$--equivariant continuous map from the Gromov boundary of $\mathcal
G(\omega_g,h,n)$ to Thurston's boundary of Teichm\"uller space, $\PML(S_g)$.

\begin{remark} We caution the reader that this notion of Cannon--Thurston map differs from that of the Kleinian
setting, where one requires that the map is a continuous extension of an equivariant map of the group into hyperbolic
space. This is {\em not} a part of our definition.\end{remark}

For the factor subgroups $G(\omega_g)$, $h^n G(\omega_g) h^{-n}$ and all their conjugates one does have
Cannon--Thurston maps in this sense.  Moreover, when a pair of these subgroups nontrivially intersect, the Cannon--Thurston maps agree on the boundary of the common subgroup.
Thus one can begin to build a ``finite approximation'' to the Cannon--Thurston map for the amalgam using these as building blocks.
Despite the existence of this approximation, we show that in general there is no such map.

\begin{theorem} \label{noct}
For all $g \geq 2$, ``generic'' $h$ as in Theorem \ref{lrthrm} and $n$ sufficiently large there does not exist a Cannon--Thurston map
$$\partial \mathcal G(\omega_g,h,n) \to \PML(S_g).$$
\end{theorem}

Here ``generic'' means that the stable lamination of $h$ lies outside a particular closed positive codimension subset
of the appropriate space of projective measured laminations---see Section \ref{examplesection} for a precise
definition. In genus $2$ for example this subset consists of just a finite set of points, and so generic means {\em not one
of a finite set of choices}, up to powers and Dehn twists in the components of $A$.

\begin{remark} In \cite{mosherproblems} Mosher considers a different boundary for $\mathcal G(\omega_g,h,n)$, viewing $\mathcal G(\omega_g,h,n)$ (more naturally) as a relatively hyperbolic group.  Then he asks about the existence of a Cannon--Thurston map to $\PML(S)$ from this boundary.  Since this boundary is a $\mathcal G(\omega_g,h,n)$--equivariant quotient of the Gromov boundary, Theorem \ref{noct} also shows that generically no Cannon--Thurston map from this boundary exists either.
\end{remark}

Whether or not there is ever a Cannon--Thurston map for a group constructed from Theorem \ref{lrthrm} remains an
open question.  However, from the proof it should be clear that such an example, if it exists, would have to be very
special.

The obstruction to extending the finite approximation to an actual Cannon--Thurston map lies in the action of parabolic mapping
classes on Thurston's boundary $\PML(S)$.  In hyperbolic geometry a parabolic isometry fixes a unique point on the
boundary of hyperbolic space and every other point is attracted to that point under iteration.  As is well known, this
is not the case for the mapping class group where there is an entire simplex of fixed points, each one ``equally
attractive'' (compare Lemma \ref{twistconvergence} below).

We mention that if we forget the transverse measures, this simplex of attractors collapses to a point. This collapse
erases the problem in building the Cannon--Thurston map, and in \cite{leiningerveech} we prove the following.
\begin{theorem} \cite{leiningerveech}
For all $g \geq 2$, all $h$ as in Theorem \ref{lrthrm} and all $n$ sufficiently large, there exists a continuous
$\mathcal G(\omega_g,h,n)$--equivariant map
$$\partial \mathcal G(\omega_g,h,n) \to \GL(S)$$
where $\GL(S)$ is the space of geodesic laminations with the Thurston topology.
\end{theorem}
This is just one part of a more general investigation into the geometry of the groups constructed in
\cite{leiningerreid}, as well as a generalization to a class of groups we call \textit{graphs of Veech groups}.\\

\noindent \textbf{Acknowledgements.}  Thanks to Richard Kent and Sergio Fenley for their interest and helpful
suggestions and conversations.

%%%%%%%%%%%%%%%%%%%%%%%%%%%%%%%%
\section{The examples revisited} \label{examplesection}
%%%%%%%%%%%%%%%%%%%%%%%%%%%%%%%%

For every $g \geq 2$, Veech constructed a genus $g$ Riemann surface together with an abelian differential $\omega_g$ by
appropriately gluing together a pair of regular $(2g+1)$--gons \cite{veech1,veech2}.  This determines a singular
Euclidean structure and an associated affine group $\Aff^+(\omega_g)$.  Locally integrating $\omega_g$ produces
preferred coordinates around all nonzero points of $\omega_g$ for which the transition functions are translations. The
derivative in these preferred coordinates defines a homomorphism
$$D: \Aff^+(\omega_g) \to \SL_2(\mathbb R)$$
The image of this map is denoted $\SL(\omega_g)$. The quotient by the center is $\PSL(\omega_g) < \PSL_2(\mathbb R)$ and is a
Fuchsian triangle group of type $(2,2g+1,\infty)$.

The homomorphism $D$ is actually an isomorphism onto $\SL(\omega_g)$ (compare \cite{leiningermult}, Section 7).  The quotient $\Aff^+(\omega_g) \to \PSL(\omega_g)$ has a central kernel of order two generated by a hyperelliptic involution we denote $\sigma$.

The commutator subgroup of $\PSL(\omega_g)$ is a free Fuchsian group with genus $g$ and one cusp.  As this is free, we can choose a lift to a free subgroup of $\Aff^+(\omega_g)$ and this is the subgroup we denote $G(\omega_g) < \Aff^+(\omega_g)$.

\begin{figure}[htb]
\begin{center}
\ \psfig{file=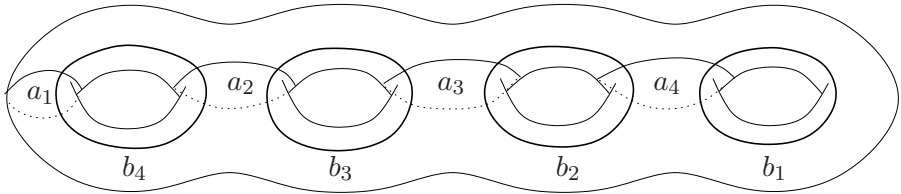,height=1truein} \caption{$A = a_1 \cup ... \cup a_4$ and $B = b_1 \cup ... \cup b_4$ in genus
$4$.} \label{aandb}
\end{center}
%%%%%%%%
  \setlength{\unitlength}{1in}
  \begin{picture}(0,0)(0,0)
    \put(-2.2,1.05){$a_1$}
    \put(-1.15,1.1){$a_2$}
    \put(-.05,1.1){$a_3$}
    \put(1.08,1.1){$a_4$}
    \put(-1.7,.65){$b_4$}
    \put(-.62,.65){$b_3$}
    \put(.57,.65){$b_2$}
    \put(1.65,.65){$b_1$}
  \end{picture}
%%%%%%%%
\end{figure}

There is a pair of cylinder decompositions for $\omega_g$, the core curves of which are multicurves
$$A = a_1 \cup ... \cup a_g \quad \mbox{ and } \quad B = b_1 \cup ... \cup b_g$$
These multicurves form a ``chain'' on $S$ as illustrated in Figure \ref{aandb}. The Dehn twists in the components of each multicurve
compose to multitwists
$$T_A = T_{a_1} \circ ... \circ T_{a_g} \quad \mbox{ and } \quad T_B = T_{b_1} \circ ... \circ T_{b_g}$$
and $T_A,T_B \in \Aff^+(\omega_g)$ with $D(T_A)$ and $D(T_B)$ parabolic.
We note that $\sigma$ is the hyperelliptic involution of $S$ leaving each of the components of $A$ and $B$ invariant.

Since $\Aff^+(\omega_g)$ has one conjugacy class of parabolic subgroups, $T_A$ and $T_B$ are conjugate in
$\Aff^+(\omega_g)$.  Similarly, the one conjugacy class of parabolic subgroups of $G(\omega_g)$ is represented by a
cyclic subgroup $G_0(\omega_g) < G(\omega_g)$ which contains a power of $T_A$.  More precisely
$$G_0(\omega_g) = \langle T_A^{2(2g+1)} \sigma \rangle.$$

Now consider $h \in \Cent(G_0(\omega_g))$.  As $h$ centralizes $T_A^{4(2g+1)} = (T_A^{2(2g+1)} \sigma)^2$ it follows
that it must leave $A$ invariant. Therefore it restricts to a mapping class on $S - A$.  As was stated in Theorem
\ref{lrthrm}, we assume that this restriction of $h$ is pseudo-Anosov.  We note that $h \in \Cent(G_0(\omega_g))$ if
and only if it leaves $A$ invariant and commutes with $\sigma$.

\begin{remark}
In \cite{leiningerreid} we erroneously stated that $G_0(\omega_g)$ was generated by a power of $T_A$, instead of by
$T_A^{2(2g+1)}\sigma$.  The reason that the generator is $T_A^{2(2g+1)} \sigma$ instead of $T_A^{2(2g+1)}$ is that for
any torsion free $1$--cusped Fuchsian group of finite area, when lifted to a subgroup of $\SL_2(\mathbb R)$ the generator of the parabolic subgroup must have $\mbox{trace} = -2$,
and not $+2$---see \cite{calegaritorus} for a nice discussion.  This only affects the possible choices of $h$, requiring
that it leave $A$ fixed {\em and also} commute with $\sigma$.
\end{remark}

We now explain the meaning of ``generic'' in the statement of Theorem \ref{noct}.  The elements $h \in
\Cent(G_0(\omega_g))$ are described (up to twists in components of $A$) as pseudo-Anosov mapping classes $\hat h$ of $S
- A$ which commute with the hyperelliptic $\hat \sigma = \sigma|_{S - A}$.  Let $X \subset \PML(S - A)$ denote the
fixed points set of $\hat \sigma$ in $\PML(S - A)$.  The space $X$ is alternatively described as the lifts of
projective classes of laminations on $(S - A)/\langle \hat \sigma \rangle$, and so $X \cong S^{2g - 3}$.  Further note
that the set of all fixed points of pseudo-Anosov mapping class $\hat h$ of $S - A$ commuting with $\hat \sigma$ is dense in
$X$ since the same is true of pseudo-Anosov fixed points in $\PML((S - A)/\langle \hat \sigma \rangle)$.  We will
construct a closed, positive codimension subset $Y \subset X$, and say that $h$ is \textbf{generic} if the stable fixed
point of $\hat h$ lies outside of $Y$.

%%%%%%%%%%%%%%%%%%%%%%%%%%%%%%%%%%%%
\section{Limit sets and parabolic fixed points} \label{limitsetsection}
%%%%%%%%%%%%%%%%%%%%%%%%%%%%%%%%%%%%

The limit set for a nonelementary subgroup $G < \Aff^+(\omega)$ of the affine group of an abelian differential $\omega$
has a fairly concrete description.  For this, note that the set of projective classes of vertical foliations for
complex multiples of $\omega$, which we denote by $\PML(\omega)$ is a circle in $\PML(S)$.  In fact, the associated
Teichm\"uller disk has an ideal boundary $\partial \mathbb H_\omega$, which by the isometry $\mathbb H^2 \cong \mathbb
H_\omega$, naturally admits a projective structure $\partial \mathbb H_\omega \cong \mathbb R \mathbb P^1$.  Moreover, there
is a natural $\Aff^+(\omega)$--equivariant, piecewise projective homeomorphism $\partial \mathbb H_\omega \to
\PML(\omega)$, see \cite{klsurvey}, Theorem 2.1.  The point is that a ray in $\mathbb H_\omega$ comes from a $1$--parameter family of Teichm\"uller
deformations where the measure on the vertical foliation of some complex multiple of $\omega$ tends to zero.
We caution the reader that according to the work of Masur \cite{twoboundaries}, the map $\partial \mathbb H_\omega \to \PML(\omega)$ is \textit{not} in general the continuous extension of the embedding of $\mathbb H_q$ into Teichm\"uller space.

Because a nonelementary subgroup $G < \Aff^+(\omega)$ contains a pseudo-Anosov mapping class, it follows from the work of McCarthy and Papadopoulos \cite{mccarthypapa} that its limit set $\Lambda_G$ is the unique minimal closed $G$--invariant set.  Since $\PML(\omega)$ is a closed
$\Aff^+(\omega)$--invariant set, we see that $\Lambda_G \subset \PML(\omega)$.  By minimality of limit sets for
Fuchsian groups, $\Lambda_G$ is precisely the image of the \textit{Fuchsian limit set} $\Lambda_{D(G)} \subset
\partial \mathbb H_\omega$.  In particular when $G$ is a Veech group---that is, a lattice such as $G(\omega_g)$---the limit set is precisely $\Lambda_G = \PML(\omega)$.

We denote the foliations of $\omega_g$ associated to the cylinder decompositions discussed above by $\nu_A$ and
$\nu_B$, respectively.  The projective
classes $[\nu_A],[\nu_B] \in \PML(\omega_g)$ lie in the simplex of measures in $\PML(S)$ determined by
$A$ and $B$, respectively.  Note that $[\nu_A]$ (respectively $[\nu_B]$) is the unique fixed point in $\PML(\omega_g)$
of $T_A$ (respectively, $T_B$).

%%%%%%%%%%%%%%%%%%%%%%%%%%%
\section{The proof}
%%%%%%%%%%%%%%%%%%%%%%%%%%%%

We require the following fact, which is certainly well-known (compare \cite{FLP}, \S 6.7 for the case of a simple
closed curve).

\begin{lemma}  \label{twistconvergence}
Let $C = c_1 \cup ... \cup c_n$ be a multicurve, $T_C = T_{c_1} \circ ... \circ T_{c_n}$ the corresponding multitwist and $[\mu] \in \PML(S)$ with $i(\mu,c_j) \neq 0$ for each $1 \leq j \leq n$.  Then
$$\lim_{k \to \infty} T_C^k([\mu]) = \left[i(\mu,c_1)c_1 + ... + i(\mu,c_n) c_n \right]$$
\end{lemma}
\begin{proof}
We assume for simplicity that no two components of $C$ are isotopic.  If this is not the case, collect together
isotopic components and replace the composition of the associated Dehn twists by a power of a single twist in one of the curves---a
composition of Dehn twists in pairwise isotopic curves is isotopic to a power of a twist in one of the curves.

In what follows we use several facts from the theory of train tracks.  We
refer the reader to \cite{pennerharer} for more on the concepts we use here.

There exists a birecurrent train track $\tau$ carrying $\mu$ so that in a neighborhood of each component $c_j$, the local picture is as in Figure \ref{train} on the left (this is one of the ``standard model'' train tracks of \cite{pennerharer} for an appropriate choice of basis multicurves, depending on $C$ and $\mu$).
Near $c_j$ the weights on the branches of $\tau$ determined by $\mu$ are $x_j = i(\mu,c_j)$, $y_j$, and $z_j = x_j + y_j$ (which is the weight on the branch not labeled in the figure).
We also see that each $c_j$ is carried by $\tau$ with the corresponding weights given by $\tilde x_j = 0$, $\tilde y_j = \tilde z_j = 1$ and all other weights zero.

\begin{figure}[htb]
\begin{center}
\ \psfig{file=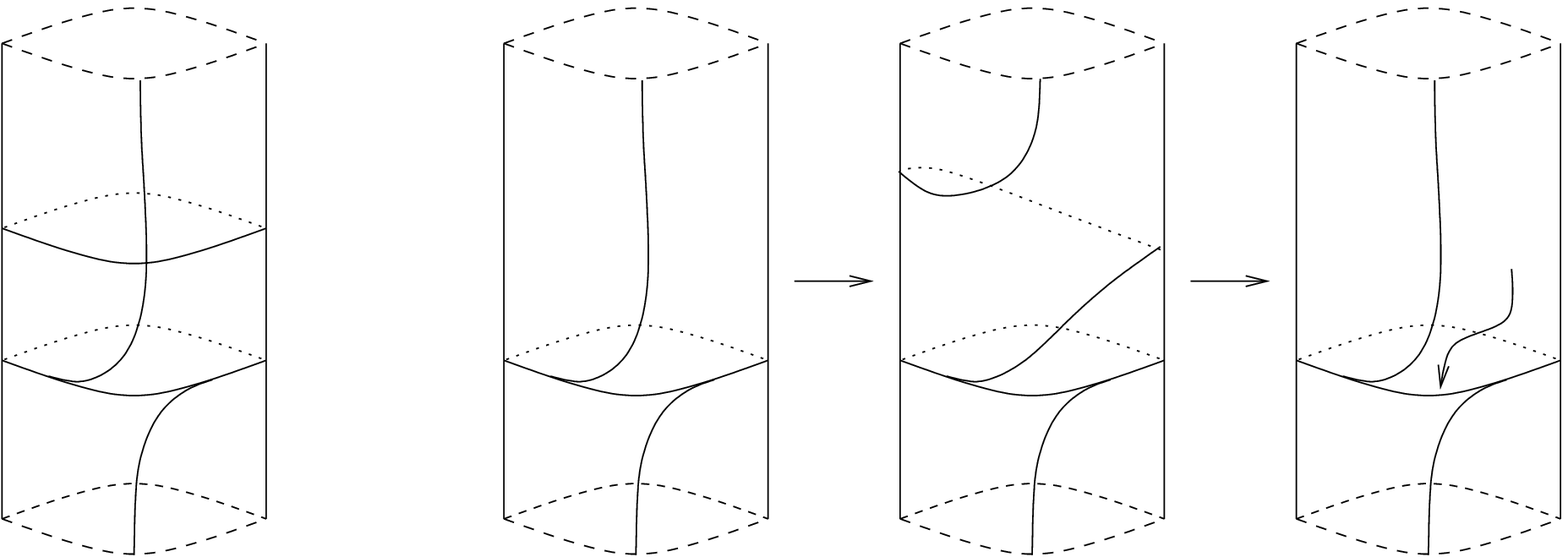,height=1.7truein}
\caption{Left: The train track near $c_j$.  Right: Dehn twisting in $C$ (first arrow) and how $\tau$ carries $T_C(\tau)$ (second arrow).}
\label{train}
\end{center}
%%%%%%%%
  \setlength{\unitlength}{1in}
  \begin{picture}(0,0)(0,0)
    \put(-2.33,1.57){$c_j$}
    \put(-2.05,1.1){$\tau$}
    \put(-.42,1.3){$y_j$}
    \put(-.58,.82){$x_j$}
    \put(-.58,2){$x_j$}
    \put(.82,1.3){$y_j$}
    \put(.62,.82){$x_j$}
    \put(.58,2){$x_j$}
    \put(2.13,1.7){$y_j + x_j$}
    \put(1.83,.82){$x_j$}
    \put(1.83,2){$x_j$}
    \put(0.1,1.7){$T_C$}
  \end{picture}
%%%%%%%%
\end{figure}

Note that $T_C(\tau)$ is carried by $\tau$, and so $T_C(\mu)$ is also carried by $\tau$.  The change in weights from $\mu$ to $T_C(\mu)$ is described on the right of Figure \ref{train}.
Observe that all weights stay the same, except each of the $y_j$ and $z_j$ which both have $x_j$ added to them.
It follows that
$$\lim_{k \to \infty} \frac{1}{k} T_C^k(\mu) = x_1 c_1 + ... + x_n c_n = i(\mu,c_1) c_1 + ... + i(\mu,c_n) c_n$$
\end{proof}

We are now prepared to prove the main theorem.
\begin{proof}[Proof of Theorem \ref{noct}]
We begin by defining the set $Y$.
First, write
\begin{equation}\label{bname}
[\nu_B] = [w_1 b_1 + ... + w_g b_g]
\end{equation}
where $w_1,...,w_g \in \mathbb R_+$ are determined by $\omega_g$---these are the heights of the cylinders as described
above. View $X \subset \PML(S - A) \subset \PML(S)$ and set
$$Y = \left\{ [\mu] \in X \, \left| \, \frac{i(\mu,b_i)}{i(\mu,b_j)} = \frac{w_i}{w_j} \mbox{ for all } 1 \leq i,j \leq g \right. \right\}$$
The set $Y$ is obtained from $X$ by imposing certain relations on intersection numbers.  Observing how the components of $B$ intersect $S - A$ (see Figure \ref{aandb}) one easily verifies that $Y$ is a positive codimension subset of $X$.

Now let $[\mu_s]$ denote the stable lamination of $h$ and assume $[\mu_s] \not \in Y$.  We suppose that there exists a Cannon--Thurston map $f:\partial \mathcal G(\omega_g,h,n) \to
\PML(S)$ and arrive at a contradiction, provided $n$ is sufficiently large.

First observe that since $\Lambda_{\mathcal G(\omega_g,h,n)}$ is the minimal closed $\mathcal G(\omega_g,h,n)$--invariant
set, $f(\partial \mathcal G(\omega_g,h,n)) = \Lambda_{\mathcal G(\omega_g,h,n)}$.

As described above, $(T_A^{2(2g+1)} \sigma) \in \mathcal G(\omega_g,h,n)$ and hence all its $\mathcal
G(\omega_g,h,n)$--conjugates also lie in $\mathcal G(\omega_g,h,n)$.  In particular, $(T_B^{2(2g+1)} \sigma) \in \mathcal
G(\omega_g,h,n)$.  Since this is the fundamental group of a closed surface, $(T_B^{2(2g+1)} \sigma)$ has exactly two
fixed points $x^\pm \in
\partial \mathcal G(\omega_g,h,n)$: $x^+$ the attracting fixed point and $x^-$ the repelling fixed point.

Since $G(\omega_g) < \mathcal G(\omega_g,h,n)$, we have
$$\PML(\omega_g) = \Lambda_{G(\omega_g)} \subset \Lambda_{\mathcal G(\omega_g,h,n)} \subset \PML(S).$$
\begin{claim}
$f(x^\pm) = [\nu_B]$.
\end{claim}
\begin{proof}
By considering the action of the parabolic $D(T_B^{2(2g+1)} \sigma)$ on $\partial \mathbb H_{\omega_g}$ and applying
the $G(\omega_g)$--equivariant homeomorphism $\partial \mathbb H_{\omega_g} \to \PML(\omega_g)$, iterating $(T_B^{2(2g+1)} \sigma)$ on $[\nu_A] \in \PML(\omega_g)$ we see
$$\lim_{k \to \pm \infty} (T_B^{2(2g+1)} \sigma)^k([\nu_A]) = [\nu_B]$$
If we let $y \in \partial \mathcal G(\omega_g,h,n)$ be an element with $f(y) = [\nu_A]$, then by continuity and $\mathcal
G(\omega_g,h,n)$--equivariance of $f$ we obtain
$$\begin{array}{rclcl}[\nu_B] & = & \displaystyle{\lim_{k \to \pm \infty} (T_B^{2(2g+1)} \sigma)^k(f(y))} & = & \displaystyle{\lim_{k \to \pm \infty} f((T_B^{2(2g+1)} \sigma)^k(y))}\\\\
 & = & \displaystyle{f(\lim_{k \to \pm \infty} (T_A^{2(2g+1)} \sigma)^k(y))} & = & f(x^\pm) \end{array}$$ proving the claim.
\end{proof}

Since $h^n G(\omega_g) h^{-n} < \mathcal G(\omega_g,h,n)$ and
$$h^n(\Lambda_{G(\omega_g)}) = \Lambda_{h^n G(\omega_g) h^{-n}} \subset \Lambda_{\mathcal G(\omega_g,h,n)}$$ it follows that
$h^n([\nu_B]) \in \Lambda_{\mathcal G(\omega_g,h,n)} = f(\partial \mathcal G(\omega_g,h,n))$.  Therefore, there exists $z
\in
\partial \mathcal G(\omega_g,h,n)$ with $f(z) = h^n([\nu_B])$, and so again by the continuity and $\mathcal
G(\omega_g,h,n)$--equivariance of $f$ as well as the claim above
$$\begin{array}{rclcl} [\nu_B] \, = \, f(x^+) & = & \displaystyle{f(\lim_{k \to \infty} (T_B^{2(2g+1)} \sigma)^k
(z))} & = & \displaystyle{\lim_{k \to \infty} (T_B^{2(2g+1)} \sigma)^k(f(z))}\\\\
& = & \displaystyle{\lim_{k \to \infty} (T_B^{2(2g+1)}\sigma)^k(h^n([\nu_B]))} & = & \displaystyle{\lim_{k \to \infty} T_B^{2k(2g+1)} (h^n([\nu_B]))} \end{array}$$
Where the last equality follows from the fact that $\sigma$ commutes with both $h$ and $T_B$, and $\sigma$ fixes $[\nu_B]$.
Appealing to Lemma \ref{twistconvergence} this implies
$$[\nu_B] = \left[ i(h^n(\nu_B),b_1) b_1 + ... + i(h^n(\nu_B),b_g) b_g \right].$$
Combining this with (\ref{bname}) we see that for all $1 \leq i,j \leq g$
\begin{equation} \label{ratio1gen}
\frac{i(h^n(\nu_B),b_i)}{i(h^n(\nu_B),b_j)} = \frac{w_i}{w_j}.
\end{equation}

Since $\mu_s$ fills $S - A$, we have $i(\nu_B,\mu_s) \neq 0$, and so by Theorem A1 of \cite{ivanov}
$$\lim_{n \to \infty} h^n([\nu_B]) = [\mu_s]$$
Therefore, for every $1 \leq i,j \leq g$
\begin{equation} \label{bslimitgen}
\lim_{n \to \infty} \frac{i(h^n(\nu_B),b_i)}{i(h^n(\nu_B),b_j)} = \frac{i(\mu_s,b_i)}{i(\mu_s,b_j)}.
\end{equation}

Since $[\mu_s] \not \in Y$, there exists indices $1 \leq i,j \leq g$ so that
$$\frac{i(\mu_s,b_i)}{i(\mu_s,b_j)} \neq \frac{w_i}{w_j}$$
Thus by (\ref{bslimitgen}), for all sufficiently large $n$
$$\frac{i(h^n(\nu_B),b_i)}{i(h^n(\nu_B),b_j)} \neq \frac{w_i}{w_j}$$
which contradicts (\ref{ratio1gen}).  This completes the proof.\\
\end{proof}

\bibliographystyle{plain}
\bibliography{multitube}

\bigskip

\noindent
Department of Mathematics, University of Illinois at Urbana--Champaign, 1409 W. Green St., Urbana, IL, 61801.\\
\texttt{clein@math.uiuc.edu}

\end{document}